\documentclass{ifacconf}

\usepackage{graphicx}      
\usepackage{natbib}        
\usepackage{amsmath} 
\usepackage{amssymb}  
\usepackage{caption}
\usepackage{subcaption}
\usepackage{xcolor}
\usepackage{tikz}
\usepackage{pgfplots}
\usepackage{lipsum}
\usepackage{comment}

\newcommand{\mcl}[1]{\mathcal{ #1}}
\newcommand{\mbf}[1]{\mathbf{ #1}}
\newcommand{\smaller}{\fontsize{9pt}{10pt}}
\def\QEDclosed{\mbox{\rule[0pt]{1.3ex}{1.3ex}}} 

\newcommand{\mylen}{3.1pt}

\newcommand{\ip}[2]{\langle{#1},{#2}\rangle}
\newcommand{\vmatwo}[2]{\begin{bmatrix}
		#1 \\
		#2
\end{bmatrix}}

\renewcommand{\th}{\ensuremath{\theta}}

\newcommand{\mc}[1]{\ensuremath{\mathcal{#1}}}

\newcommand{\myint}{\int_{a}^{b}}
\newcommand{\myinta}[1]{\int_{a}^{#1}}
\newcommand{\myintb}[1]{\int_{#1}^{b}}

\newcommand{\ltwo}{\ensuremath{L_2}}
\newcommand{\rl}{\ensuremath{\mathbb{R}}}

\newcommand{\R}{\mathbb{R}}

\renewcommand{\S}{\mathbb{S}}
\newcommand{\opPzero}{\ensuremath{\mc{P}}}
\newcommand{\opPgen}{\ensuremath{\mc{P}_{\{P,Q_1,Q_2,S,R_1,R_2\}}}}

\newcommand{\opPpq}[4]{\ensuremath{\mc{P}_{\tiny\{#1,#2_1,#2_2,#4,#3_1,#3_2\}}}}

\newcommand{\opPselfadj}[4]{\ensuremath{\mc{P}_{\tiny\{#1,#2,#2^{\top},#4,#3_1,#3_2\}}}}

\allowdisplaybreaks
\begin{document}
\begin{frontmatter}

\title{Representation and Stability Analysis of PDE-ODE Coupled Systems} 

\author[First]{Amritam Das} 
\author[Second]{Sachin Shivakumar} 
\author[First]{Siep Weiland} 
\author[Second]{Matthew Peet}

\address[First]{Department of Electrical Engineering, Eindhoven University of Technology, 5600 MB, The Netherlands (e-mail: \{am.das,s.weiland\}@tue.nl).}
\address[Second]{School for Engineering of Matter, Transport and Energy, Arizona State University, 
   Tempe, AZ 85292 USA (e-mail: \{sshivak8, mpeet\}@asu.edu)}

\begin{abstract}                
In this work, we present a scalable Linear Matrix Inequality (LMI) based framework to verify the stability of a set of linear Partial Differential Equations (PDEs) in one spatial dimension coupled with a set of Ordinary Differential Equations (ODEs) via input-output based interconnection. Our approach extends the newly developed state space representation and stability analysis of coupled PDEs that allows parametrizing the Lyapunov function on $L_2$ with multipliers and integral operators using polynomial kernels of semi-separable class. In particular, under arbitrary well-posed boundary conditions, we define the linear operator inequalities on $\mathbb{R}^n \times L_2$ and cast the stability condition as a feasibility problem constrained by LMIs. In this framework, no discretization or approximation is required to verify the stability conditions of PDE-ODE coupled systems.  The developed algorithm has been implemented in MATLAB where the stability of example PDE-ODE coupled systems are verified.
\end{abstract}

\begin{keyword}
Partial Differential Equations, Ordinary Differential equations, Stability, Lyapunov Function, Linear Matrix Inequality, Sum of Squares.
\end{keyword}

\end{frontmatter}
\section{Introduction}
\begin{figure}[h!]
	\begin{center}
		\includegraphics[width=4cm]{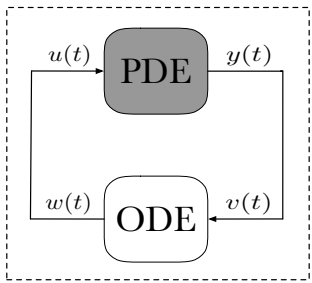}   
		\caption{Interconnection between PDE-ODE models.} 
		\label{fig:Feedback_interconnection}
	\end{center}
\end{figure}

Coupled Ordinary Differential Equations(ODEs)-Partial Differential Equations(PDEs) are useful to model systems that involve states varying over time as well as states varying over both space and time. Interaction among systems that are governed by ODEs and PDEs appear in multi-scale modeling of fluid flows, for instance in \cite{quarteroni2003analysis} and \cite{stinner2014global}, the flow structures are the function of space and time; however, they are driven by global flow pattern that is only a function of time. They also frequently appear in thermo-mechanical systems (e.g. \cite{tang}) and chemical reactions (e.g. \cite{christofides2012nonlinear}). A large class of these applications involves linear ODEs and PDEs that are mutually coupled via inputs and outputs as illustrated in Fig. \ref{fig:Feedback_interconnection}. 

In this paper, the linear ODEs are of the following form:
\begin{equation}
\label{intro1}
\begin{split}
&\dot{x}(t) = A x(t) + B v(t),\\ 
&w(t) = C x(t).
\end{split}
\end{equation}
For $s \in [a, b]$, the linear one dimensional PDEs are described in the following form:
\begin{align}
\label{intro2}
&\dot{z}(s,t) = A_0(s)z(s,t) + A_1(s)\frac{\partial z(s,t)}{\partial s} \\
&\qquad\qquad\qquad+A_2(s)\frac{\partial^2 z(s,t)}{\partial s^2} + B_1(s)u(t), \nonumber\\
&y(t)= C_1 z(t) + \int_{a}^b \left(C_a(s)z(s,t)+C_b(s)\frac{\partial z(s,t)}{\partial s}\right)\text{d}s, \nonumber\\
&B_c \text{col}\Big(z(a,t), z(b,t), \frac{\partial z(s,t)}{\partial s}\mid_{s = a}, \frac{\partial z(s,t)}{\partial s}\mid_{s=b}\Big) = 0.  \nonumber
\end{align}
Here, $x(t)\in\R^{n_{\text{o}}}$ and $z(s,t)\in\R^{n_{\text{p}}}$ are the states of the ODEs and the PDEs respectively. The models \eqref{intro1}-\eqref{intro2} are connected by the signals $\text{col}(v(t), w(t)) \in \R^{m_{\text{o}}}\times \R^{p_{\text{o}}}$ and $\text{col}(u(t), y(t)) \in \R^{m_{\text{p}}}\times \R^{p_{\text{p}}}$. In \eqref{intro1} $A, B, C$ are constant matrices of suitable dimensions. In \eqref{intro2}, $A_0(s), A_1(s),A_2(s), B_1(s), C_a(s), C_b(s)$ are matrix valued functions of appropriate dimensions and $C_1$ is a constant matrix. $B_c \in \mathbb{R}^{2n_{\text{p}} \times 4n_{\text{p}}}$ is required to have full row rank. 

The interconnection between \eqref{intro1}-\eqref{intro2} is described by a matrix $L \in \mathbb{R}^{{m_{\text{o}}+m_{\text{p}}} \times {p_{\text{o}}+p_{\text{p}}}}$ that has full row rank and relates the interconnection signals with the help of the following algebraic equation:
\begin{equation}
\left[ \begin{array}{c}
v(t)   \\
u(t)  
\end{array}\right]=
L\left[ \begin{array}{c}
w(t)   \\
y(t)  
\end{array}\right].
\end{equation}
In this paper, we develop a scalable approach to verify the stability of the coupled system in \eqref{intro1}-\eqref{intro2} that does not depend on any discretization or approximation. Regarding verifying stability, many methods used in the analysis of Distributed Parameter Systems(DPS) rely on either application-specific analytical approaches(e.g. \cite{nicaise2009stability}) or projection of the entire system onto a finite-dimensional vector space as shown in \cite{ravindran2000reduced}, \cite{das2018model} and the references therein. Once projected to a finite-dimensional subspace, the reduced order model is treated as a set of ODEs, and there are many well-studied approaches to test for stability and design controllers. However, discretization techniques for approximating an infinite-dimensional system with a finite-dimensional model are prone to numerical instabilities and other truncation errors. This results in large system matrices which often make them computationally demanding. Moreover, these finite-dimensional approximations of the DPS do not necessarily represent the behavior of the original system to a required level of accuracy, thus the properties, such as stability, robustness are not exact if the finite dimensional approximation replaces the original model. Another approach frequently used in stabilization and control of coupled PDE-ODE systems is the backstepping method, e.g. \cite{krstic2008backstepping}, \cite{krstic2009compensating} and \cite{susto2010control}. The backstepping technique does not require any discretization or approximation of the model; however, it does not allow us to utilize Lyapunov theory to verify stability. Also, these methods are neither generic nor scalable for a general class of PDE-ODE coupled system. 

On the other hand, using the extension of Lyapunov theory in infinite-dimensional space by \cite{datko1970extending}, many attempts were made to use Sum-of-Squares (SOS) optimization methods for constructing Lyapunov functions for PDEs. Some of the notable works include but not limited to \cite{Papa}, \cite{fridman2009exponential} and \cite{valmorbida2016stability}. However, these works cannot deal with a general class of coupled PDEs that involve reaction, transport, and diffusion terms. From our lab, numerous problems related to the analysis of stability and robustness and controller-observer synthesis are addressed in \cite{meyer2015stability}, \cite{Gahlawat} etc. Recently, \cite{peet2018new}, \cite{Sachin} have developed a novel SOS framework for representation, stability and robustness analysis for coupled PDEs with arbitrary boundary condition. These results provide a unifying framework and a computationally scalable method for the analysis of a large class of PDEs.

The contributions of this paper involve extending the representation and the stability analysis of coupled PDEs to the case of linear PDE-ODE coupled systems. To this end, we construct a quadratic Lyapunov function of the form $V({x},{z})=\Big\langle{\vmatwo{x}{z}},{\opPselfadj{P}{Q}{R}{S}\vmatwo{x}{z}}\Big\rangle$, where $\opPselfadj{P}{Q}{R}{S}$ 
is a class of self-adjoint linear operator on $\R^{n_{\text{o}}}\times\ltwo$. The criterion for the stability of the PDE-ODE coupled system is to search for a positive \opPgen that restricts the time derivative of the Lyapunov function to be negative. We parametrize the multipliers and integral operators of \opPgen with polynomial kernels of semi-separable class as in \eqref{intro_operator}. We show that the positivity constraints on these operators can be re-formulated as scalable Linear Matrix Inequalities(LMIs) that establish the stability condition for PDE-ODE coupled system. 
{\small
	\begin{align}
	\label{intro_operator}
	&\left(\opPselfadj{P}{Q}{R}{S}\vmatwo{x}{z}\right)(s):= \nonumber\\
	& {\scriptsize\vmatwo{Px + \myint Q(s) z(s)ds}{Q^{\top}(s)x + S(s)z(s)+\myinta{s}R_1(s,\eta)z(\eta)\text{d}\eta+\myintb{s}R_2(s,\eta)z(\eta)\text{d}\eta}}
	\end{align}}
The organization of the paper is as follows. Section \ref{section:preliminaries} provides a preliminary discussion on the notations, a model of the PDE-ODE coupled systems and the Lyapunov stability theory. In Section \ref{section:operator_def}, we define the class of operators \opPgen and derive its positivity condition as well as the composition and adjoint operations of these operators. In Section \ref{section:fundamental_identities}, the PDE-ODE coupled model has been represented using the same class of operators. In Section \ref{section:stability}, the equivalent LMI conditions have been derived to verify the stability of PDE-ODE coupled system. In Section \ref{section:results}, we numerically illustrate the methodology to verify the stability of example PDE-ODE coupled systems in MATLAB. At last, Section \ref{section:conclusions} provides some conclusions and directions for future research.
\section{Preliminaries}
\label{section:preliminaries}
\subsection{Notation}
For convenience, we denote $\frac{\partial x}{\partial s}:= x_{s}$ and $\frac{\partial^2 x}{\partial s^2}:= x_{ss}$. We use $\mathbb{Z}$ to denote a field of integers. We use $\mathbb{S}^m\subset \rl^{m\times m}$ to denote the symmetric matrices. We define the space of square integrable $\rl^m$-valued functions on $X$ as $\ltwo^m(X)$. $\ltwo^m(X)$ is equipped with the inner product $\langle x,y \rangle_{\ltwo} = \myint x^{\top}(s) y(s) \text{d}s$. For denoting the inner product in $\mathbb{R}^{n}$, we use $\ip{x}{y}_{\mathbb{R}^{n}}=x^{\top}y$. The Sobolov space is defined by $W^{q,n}(X):=\{x\in \ltwo^n(X) \mid \frac{\partial^k x}{\partial s^k}\in\ltwo^n(X) \text{ for all }k\le q \}$. For an inner product space $X$, operator $\mc{P}:X\to X$ is called positive, if for all $x\in X$, we have $\ip{x}{\mc{P}x}_X\ge 0$. We use $\mc{P}\succcurlyeq 0$ to indicate that $\mcl P$ is a positive operator. We say that $\mcl P:X\to X$ is coercive if there exists some $\epsilon>0$ such that $\ip{x}{\opPzero x}_X\ge \epsilon \Vert x \Vert_X^2$ for all $x \in X$.
\subsection{PDE-ODE Coupled System}
In this section, we elaborate on the description of the coupled linear PDE-ODE system.
The system is represented as an interconnection between two systems; one of them is governed by a set of partial differential equations (PDEs) while the other is governed by a set of ordinary differential equations (ODEs). In the following subsections we briefly describe the corresponding PDE and ODE model. 
\subsubsection{A. ODE Model:}
The ODE model has the following state space representation:
\begin{equation}
    \label{eq:model_ODE}
    \begin{split}
    \dot{x}(t) &= A x(t) + B v(t),\\ 
    w(t) &= C x(t),
    \end{split}
\end{equation}
where $x(t) \in \mathbb{R}^{n_{\text{o}}}$ are the state variables, $v(t) \in \mathbb{R}^{m_{\text{o}}}$, $w(t) \in \mathbb{R}^{p_{\text{o}}}$ are the interconnection signals. $A, B, C$ are matrices of suitable dimensions.  
\subsubsection{B. PDE Model:}
An operator based notation can be utilized to represent the PDE model in the following abstract state space form:
\begin{equation}
    \label{eq:model_PDE}
    \begin{split}
    \dot{z}(t) &= \mathcal{A} z(t) + \mathcal{B} u(t),\\ 
    y(t) &= \mathcal{C} z(t).
        \end{split}
\end{equation}
For functions $\mathbf{z}: [a, b] \times \mathbb{R}^{+} \rightarrow \mathbb{R}^{n_{\text{p}}}$, $\mathbf{u}: \mathbb{R}^{+} \rightarrow \mathbb{R}^{m_{\text{p}}}$, $\mathbf{y}: \mathbb{R}^{+} \rightarrow \mathbb{R}^{p_{\text{p}}}$, the system operators $\mathcal{A}:L_2^{n_{\text{p}}}[a, b]\supset D_{\mathcal{A}}\to L_2^{n_{\text{p}}}[a, b], ~\mathcal{B}:\rl^{m_\text{p}}\to L_2^{n_{\text{p}}}[a, b]$ and $~\mathcal{C}:L_2^{n_{\text{p}}}[a, b]\to\rl^{p_\text{p}}$ are specified by the following definitions:
\begin{align}
    \label{PDE_dynamics}
    (\mathcal{A}\mathbf{z})(s) &:=  A_0(s)\mathbf{z}(s) + A_1(s)\mathbf{z}_{s}+ A_2(s) \mathbf{z}_{ss},\notag\\
    (\mc{B}\mathbf{u})(s) &:= B_1(s)\mathbf{u},\nonumber\\
    \mc{C}\mathbf{z} &:=  C_1 \mathbf{z}_b + \int\limits_{a}^{b} \Bigg(C_a(s) \mathbf{z} + C_b(s) \mathbf{z}_s\Bigg) \text{d}s,
\end{align}
with its domain 
\begin{equation}
\label{PDE_boundary}
\begin{split}
    D_{\mathcal{A}}:=&\{\mathbf{z}\in W^{2,2}: B_c \mathbf{z}_b = 0 \},\\
    \mathbf{z}_b:=&\text{col}\big(\mathbf{z}(a), \mathbf{z}(b), \mathbf{z}_{s}(a), \mathbf{z}_{s}(b)\big).
    \end{split}
\end{equation}
In \eqref{PDE_dynamics}, $A_0(s), A_1(s), A_2(s) \in \mathbb{R}^{n_{\text{p}} \times n_{\text{p}}}$ are matrix valued functions. $B_1(s) \in \mathbb{R}^{n_{\text{p}} \times m_{\text{p}}}$, $C_a(s), C_b(s) \in \mathbb{R}^{m_{\text{p}} \times n_{\text{p}}}$ are also matrix valued functions. $C_1 \in \mathbb{R}^{m_{\text{p}}\times 4n_{\text{p}}}$ is a constant matrix. $B_c \in \mathbb{R}^{2n_{\text{p}}\times 4n_{\text{p}}}$ has full row rank.
\subsubsection{C. PDE-ODE Interconnection:}
A regular interconnection between the ODE model in \eqref{eq:model_ODE} and the PDE model in \eqref{eq:model_PDE} is specified by the following algebraic constraint on the variable pairs $(v, w)$ and $(u, y)$:
\begin{equation}
\label{interconnection_matrix}
   \left[ \begin{array}{c}
        v(t)   \\
        \hline
        u(t)  
    \end{array}\right]=
\underbrace{\left[
\begin{array}{c|c}
L_1 & L_2 \\
\hline
L_3 & L_4
\end{array}
\right]}_{L}\left[ \begin{array}{c}
        w(t)   \\
        \hline
        y(t) 
    \end{array}\right].
\end{equation}
The matrix $L \in \mathbb{R}^{{m_{\text{o}}+m_{\text{p}}} \times {p_{\text{o}}+p_{\text{p}}}}$ has full row rank and is partitioned to relate interconnection signals of ODEs and PDEs individually. Substituting \eqref{interconnection_matrix} in \eqref{eq:model_ODE} and \eqref{eq:model_PDE}, we obtain the following autonomous state-space model.
\begin{equation}
    \label{combined-system}
    \begin{bmatrix}
\dot{\mathbf{x}}(t)\\
\dot{\mathbf{z}}(t)
    \end{bmatrix} = \underbrace{\begin{bmatrix}
    A+BL_1C & B L_2\mathcal{C}\\
    \mathcal{B}L_3C & \mathcal{A}+ \mathcal{B}L_4\mathcal{C}
    \end{bmatrix}}_{\bar{\mbf{A}}}\begin{bmatrix}
\mathbf{x}(t)\\
\mathbf{z}(t)
    \end{bmatrix},
\end{equation}
with $\bar{\mbf{A}}: \mathbb{R}^{n_{\text{o}}} \times L_2^{n_{\text{p}}}[a, b] \supset \mathbb{R}^{n_{\text{o}}}\times D_{\mathcal{A}}\rightarrow \mathbb{R}^{n_{\text{o}}} \times L_2^{n_{\text{p}}}[a, b]$. We define $X := \mathbb{R}^{n_{\text{o}}} \times L_2^{n_{\text{p}}}[a, b]$ with inner product
\begin{align}
\Bigg\langle\begin{bmatrix}
x_1\\
x_2(s)
\end{bmatrix}, \begin{bmatrix}
z_1\\
z_2(s)
\end{bmatrix}\Bigg\rangle_{X}:= x_1^{\top} z_1 + \big\langle x_2, z_2\big\rangle_{L_2^{n_{\text{p}}}[a, b]},
\end{align} 
and norm
\begin{align}
\Bigg\lvert\Bigg\lvert\begin{bmatrix}
    x_1\\
    x_2(s)
    \end{bmatrix}\Bigg\rvert\Bigg\rvert_{X}^{2} = \mid\mid x_1\mid\mid_{\mathbb{R}^{n_{\text{o}}}}^2 + \mid\mid x_2(s)\mid\mid_{L_2^{n_{\text{p}}}[a, b]}^2.
\end{align}
\subsection{Lyapunov Stability Theorem}
The necessary and sufficient condition for exponential stability of a linear infinite dimensional system has been provided in \cite{Curtain:1995:IIL:207416}. The same can be extended as Theorem \ref{th:stability_condition} that provides a condition to verify the exponential stability of PDE-ODE coupled system. 
\begin{thm}
\label{th:stability_condition}
    Suppose $\bar{\mbf{A}}: X \supset \mathbb{R}^{n_{\text{o}}}\times D_{\mathcal{A}}\rightarrow X$ generates a strongly continuous semigroup and there exists $\alpha, \beta, \gamma>0$ and $\opPzero: X\to X$ such that 
    \[
    \alpha\Bigg\lvert\Bigg\lvert\begin{bmatrix}
    x\\
    z
    \end{bmatrix}\Bigg\rvert\Bigg\rvert_{X}^{2} \le \Bigg \langle\begin{bmatrix}
    x\\
    z
    \end{bmatrix}, {\opPzero \begin{bmatrix}
    x\\
    z
    \end{bmatrix}}\Bigg \rangle_{X}\le\beta\Bigg\lvert\Bigg\lvert\begin{bmatrix}
    x\\
    z
    \end{bmatrix}\Bigg\rvert\Bigg\rvert_{X}^{2},\] and 
 
    \begin{align}
      \Bigg \langle  \begin{bmatrix}
    x\\
    z
    \end{bmatrix},{\opPzero \bar{\mbf{A}}\begin{bmatrix}
    x\\
    z
    \end{bmatrix}}\Bigg \rangle_{X} + \Bigg \langle \bar{\mbf{A}}\begin{bmatrix}
    x\\
    z
    \end{bmatrix},{\opPzero \begin{bmatrix}
    x\\
    z
    \end{bmatrix}}\Bigg \rangle_{X}\le -\gamma\Bigg\lvert\Bigg\lvert\begin{bmatrix}
    x\\
    z
    \end{bmatrix}\Bigg\rvert\Bigg\rvert_{X}^{2},
    \end{align}
for all $x \in \mathbb{R}^{n_{\text{o}}}$, $z\in D_{\mathcal{A}}$. Then \eqref{combined-system} is exponentially stable.
\end{thm}

\section{Parametrization of the Class of Operator \opPgen}
\label{section:operator_def}
\begin{defn}
\label{def_operator}
For a matrix $P \in \mathbb{R}^{n_{\text{o}}\times n_{\text{o}}}$, and bounded polynomial functions $Q_1: [a, b] \rightarrow \mathbb{R}^{n_{\text{o}}\times n_{\text{p}}}$, $Q_2: [a, b] \rightarrow \mathbb{R}^{n_{\text{p}}\times n_{\text{o}}}$, $S: [a, b] \rightarrow \mathbb{R}^{n_{\text{p}}\times n_{\text{p}}}$, and $R_1, R_2: [a, b]\times [a, b] \rightarrow \mathbb{R}^{n_{\text{p}}\times n_{\text{p}}}$ with $x \in \mathbb{R}^{n_{\text{o}}}$ and $z\in L_2^{n_{\text{p}}}[a, b]$, the class of operators \opPgen are defined as
{
\begin{align}
\label{mathcalP}
    &\left(\opPgen\vmatwo{x}{z}\right)(s):= \nonumber\\
    & {\scriptsize\vmatwo{Px + \myint Q_1(s) z(s)ds}{Q_2(s)x + S(s)z(s)+\myinta{s}R_1(s,\eta)z(\eta)\text{d}\eta+\myintb{s}R_2(s,\eta)z(\eta)\text{d}\eta}}.
\end{align}}
\end{defn}
In the remainder of this section, we investigate three properties of the class of operator \opPgen. First, we provide a sufficient condition on the positivity of the operator. Subsequently, an efficient method to compute composition and adjoint of the operators of the form \opPgen has been provided. The composition and the adjoint operation on the operators will allow us to obtain a scalable represention of the PDE-ODE coupled system. With the help of the positivity condition on the operator, the Lyapunov stability conditions can be reformulated as LMIs.

\subsection{Positivity of \opPselfadj{P}{Q}{R}{S} Operator}
\begin{thm}\label{th:positivity}
    For any functions $Z(s)\in\R^{d_1\times n}$, $Z(s,\eta)\in\mathbb{R}^{d_2\times n}$ $\forall s, \eta \in [a, b]$, suppose there exists a matrix 
    $T\succcurlyeq0$ such that
    \setlength{\abovedisplayskip}{\mylen}
    \setlength{\belowdisplayskip}{\mylen}
   \begin{align}
    P &= T_{11},\nonumber\\
    Q(s) &= T_{12}Z(s)+\myintb{s} T_{13}Z(\eta,s)\text{d}\eta+\myinta{s} T_{14}Z(\eta,s)\text{d}\eta, \nonumber\\
    S(s) &= Z^{\top}(s) T_{22} Z(s),\nonumber\\
    R_1(s,\eta) &=Z^{\top}(s)T_{23}Z(s,\eta)+Z^{\top}(\eta,s)T_{42}Z^{\top}(\eta)\notag\\
    &\hspace{-1.1cm}+\myintb{s}Z^{\top}(\theta,s)T_{33}Z(\theta,\eta)\text{d}\theta+\int_{\eta}^{s}Z^{\top}(\theta,s)T_{43}Z(\theta,\eta)\text{d}\theta\nonumber\\
    &\hspace{-1cm}+\myinta{\eta}Z^{\top}(\theta,s)T_{44}Z(\theta,\eta)\text{d}\theta,\nonumber\\
    R_2(s,\eta) &=Z^{\top}(s)T_{32}Z(s,\eta)+Z^{\top}(\eta,s)T_{24}Z^{\top}(\eta)\notag\\
    &\hspace{-1.1cm}+\myintb{\eta}Z^{\top}(\theta,s)T_{33}Z(\theta,\eta)\text{d}\theta+\int_{s}^{\eta}Z^{\top}(\theta,s)T_{34}Z(\theta,\eta)\text{d}\theta\nonumber\\
    &\hspace{-1cm}+\myinta{s}Z^{\top}(\theta,s)T_{44}Z(\theta,\eta)\text{d}\theta,
    \label{eq:TH}
   \end{align}
with
   \begin{align*}
   T= \begin{bmatrix}
   T_{11} & T_{12} & T_{13} & T_{14}\\
   T_{21} & T_{22} & T_{23} & T_{24}\\
   T_{31} & T_{32} & T_{33} & T_{34}\\
   T_{41} & T_{42} & T_{43} & T_{44}
   \end{bmatrix}.
   \end{align*}\\
   Then the operator \opPselfadj{P}{Q}{R}{S} ~as defined in \eqref{mathcalP} is positive, i.e. $\opPselfadj{P}{Q}{R}{S}\succcurlyeq 0$.
   \end{thm}
   \begin{pf}
   \setlength{\abovedisplayskip}{\mylen}
    \setlength{\belowdisplayskip}{\mylen}
   Since $T\succcurlyeq 0$, we can define a square root of $T$ as $U = \begin{bmatrix} U_1&U_2&U_3&U_4\end{bmatrix}$.
   \begin{align*}\label{eq:Teq}
   T = \begin{bmatrix}
   T_{11} & T_{12} & T_{13} & T_{14}\\
   T_{21} & T_{22} & T_{23} & T_{24}\\
   T_{31} & T_{32} & T_{33} & T_{34}\\
   T_{41} & T_{42} & T_{43} & T_{44}
   \end{bmatrix} &= U^{\top} U\\
   &=\begin{bmatrix}
   U_1^{\top}U_1 & U_1^{\top}U_2 & U_1^{\top}U_3 &U_1^{\top}U_4\\
   *^{\top} & U_2^{\top}U_2 & U_2^{\top}U_3 &U_2^{\top}U_4\\
   *^{\top} & *^{\top} & U_3^{\top}U_3 &U_3^{\top}U_4\\
   *^{\top} & *^{\top} & *^{\top} &U_4^{\top}U_4
   \end{bmatrix}
   \end{align*}
For $x_1 \in \mathbb{R}^{m}$ and $x_2 \in L_2^n[a, b]$, let us define $v(s) := U_1x_1 + (\Psi x_2)(s)$, where
\begin{align*}
(\Psi x_2)(s)&=U_2 Z(s) + \myinta{s} U_3 Z(s,\eta) x_2(\eta) \text{d}\eta\\
&\quad+\myintb{s} U_4 Z(s,\eta) x_2(\eta) \text{d}\eta.
\end{align*}
 Then
   \begin{align*}
       &\Big\langle{v(s)},{v(s)}\Big\rangle_{\R^{m}\times L_2^n[a, b]}\\
       &= \ip{U_1x_1}{U_1x_1}_{\mathbb{R}^{m}}+\ip{U_1x_1}{(\Psi x_2)(s)}_{\mathbb{R}^{m}}\\
   &~~~+\ip{(\Psi x_2)(s)}{U_1x_1}_{L_2^n[a, b]}+\ip{(\Psi x_2)(s)}{(\Psi x_2)(s)}_{L_2^n[a, b]}\\
   &\quad=\Big\langle{\vmatwo{x_1}{x_2}}{\opPselfadj{P}{Q}{R}{S}, \vmatwo{x_1}{x_2}}\Big\rangle_{\R^{m}\times L_2^n[a, b]}\ge 0.
   \end{align*}
Hence $\opPselfadj{P}{Q}{R}{S} \succcurlyeq 0$.
\QEDclosed
\end{pf}
\subsection{Composition of The Operators}
In this subsection, we show that the composition of two operators of the class \opPgen also takes the same structure and provide an efficient method to calculate it.
\begin{lem}\label{lem:composition}
For any matrices $A,P\in\R^{m\times m}$ and bounded functions $B_1, Q_1: [a, b]\to\R^{m\times n}$, $B_2, Q_2: [a, b]\to\R^{n\times m}$, $D,S:[a, b]\to\R^{n \times n}$, $C_i,R_i:[a, b]\times [a, b] \to \R^{n\times n}$ with $i \in \{1,2\}$, the following identity holds 
\begin{align}
    \opPpq{A}{B}{C}{D}&\opPpq{P}{Q}{R}{S} \nonumber\\
    &\hspace{-1cm}= \mc{P}_{\{\hat{P},\hat{Q}_1,\hat{Q}_2,\hat{S},\hat{R}_1,\hat{R}_2\}},
\end{align}
where
\begin{align}
\label{compostion_matrix}
    \hat{P} &= AP + \myint B_1(s)Q_2(s)\text{d}s,\nonumber\\
    \hat{Q}_1(s) &= AQ_1(s) + B_1(s)S(s)+\myintb{s}B_1(\eta)R_1(\eta,s)\text{d}\eta\nonumber\\
    &\qquad+\myinta{s}B_1(\eta)R_2(\eta,s)\text{d}\eta,\nonumber\\
    \hat{Q}_2(s) &= B_2(s)P + D(s)Q_2(s) + \myinta{s}C_1(s,\eta)Q_2(\eta)\text{d}\eta\nonumber\\
    &\qquad+\myintb{s}C_2(s,\eta)Q_2(\eta)\text{d}\eta,\nonumber\\
    \hat{S}(s) &= D(s)S(s),\nonumber\\
    \hat{R}_1(s,\eta) &=B_2(s)Q_1(\eta)+D(s)R_1(s,\eta)+C_1(s,\eta)S(\eta)\nonumber\\
    &\hspace{-0.5cm}+\myinta{\eta} C_1(s,\theta)R_2(\theta,\eta)\text{d}\theta+\int_{\eta}^{s}C_1(s,\theta)R_1(\theta,\eta)\text{d}\theta\nonumber\\
    &\hspace{-0.5cm}+\myintb{s}C_2(s,\theta)R_1(\theta,\eta)\text{d}\theta,\nonumber\\
    \hat{R}_2(s,\eta) &=B_2(s)Q_1(\eta)+D(s)R_2(s,\eta)+C_2(s,\eta)S(\eta)\nonumber\\
    &\hspace{-0.5cm}+\myinta{s} C_1(s,\theta)R_2(\theta,\eta)\text{d}\theta+\int_{s}^{\eta}C_2(s,\theta)R_2(\theta,\eta)\text{d}\theta \nonumber\\
    &\hspace{-0.5cm}+\myintb{\eta}C_2(s,\theta)R_1(\theta,\eta)\text{d}\theta.
\end{align}    
\end{lem}
\begin{pf}
The proof can be found in the Appendix C.
\QEDclosed
\end{pf}
\begin{notation}
For convenience, we say
\begin{align}
    \Big(\hat{P}, \hat{Q}_1, \hat{Q}_2, \hat{S}, \hat{R}_1, \hat{R}_2\Big) =& \Big(A, B_1, B_2, D, C_1, C_2\Big)\\
    &\times \Big(P, Q_1, Q_2, S, R_1, R_2\Big),
\end{align}
if $\hat{P}, \hat{Q}_1, \hat{Q}_2, \hat{S}, \hat{R}_1, \hat{R}_2$ are related to $A, B_1, B_2, D, C_1, C_2$ and $P, Q_1, Q_2, S, R_1, R_2$ according to \eqref{compostion_matrix} in Lemma \ref{lem:composition}.
\end{notation}
\subsection{Adjoint of The Operators}
Now, we also discuss the adjoint of an operator of the class \opPgen and it is typically denoted by $\mathcal{P}_{\{P, Q_1, Q_2, S, R_1, R_2\}}^{*}$.
\begin{lem}\label{lem:adjoint}
For any matrices $P\in\R^{m\times m}$ and  bounded functions $Q_1: [a, b]\to\R^{m\times n}$, $Q_2: [a, b]\to\R^{n\times m}$, $S:[a, b]\to\R^{n \times n}$, $R_1, R_2:[a, b]\times [a, b] \to \R^{n\times n}$,  the following identity holds for any $x \in \mathbb{R}^{m}, y \in L_2^n([a, b])$.
\begin{align}
    &\Big\langle{x},{\opPpq{P}{Q}{R}{S}y}\Big\rangle_{\R^{m}\times L_2^n[a, b]}\nonumber\\
    &=\Big\langle{\opPpq{P}{Q}{R}{S}^*x},{y}\Big\rangle_{\R^{m}\times L_2^n[a, b]},
\end{align}
where, $\opPpq{P}{Q}{R}{S}^*=\opPpq{\hat{P}}{\hat{Q}}{\hat{R}}{\hat{S}}$ with
\begin{align}
\label{adjoint_matrix}
    &\hat{P} = P^{\top}, &\hat{S}(s) = S^{\top}(s), \nonumber\\
    &\hat{Q}_1(s) = Q_2^{\top}(s), \nonumber &\hat{R}_1(s,\eta) = R_2^{\top}(\eta,s), \nonumber\\
    &\hat{Q}_2(s) = Q_1^{\top}(s),  &\hat{R}_2(s,\eta) = R_1^{\top}(\eta,s). 
\end{align}
\end{lem}
\begin{pf}
The proof can be found in the Appendix D.
\QEDclosed
\end{pf}
\begin{notation}
For convenience, we say
\begin{align}
    \Big(\hat{P}, \hat{Q}_1, \hat{Q}_2, \hat{S}, \hat{R}_1, \hat{R}_2\Big) =\ \Big(P, Q_1, Q_2, S, R_1, R_2\Big)^*,
\end{align}
if $\hat{P}, \hat{Q}_1, \hat{Q}_2, \hat{S}, \hat{R}_1, \hat{R}_2$ are related to $P, Q_1, Q_2, S, R_1, R_2$ according to \eqref{adjoint_matrix} in Lemma \ref{lem:adjoint}.
\end{notation}
\section{Fundamental Identities}
\label{section:fundamental_identities}
In this section, we utilize the class of operators defined by \opPgen to uniquely represent the system \eqref{combined-system} in terms of the variables $\mathbf{x}$ and $\mathbf{z}_{ss}$ replacing $\mathbf{z}$. 
\begin{lem}\label{lemma_xs_x}
Suppose $\mathbf{x} \in \mathbb{R}^{n_{\text{o}}}$ and $\mathbf{z} \in W^{2,2}[a, b]$ and $$B_c \text{col}(\mathbf{z}(a), \mathbf{z}(b), \mathbf{z}_s(a), \mathbf{z}_s(b)) = 0,$$ with $B_c\in \R^{2n_{\text{p}}\times 4n_{\text{p}}}$ has full row rank. Then 
\begin{equation}
\label{basic_identity_x}
      \begin{bmatrix}
      \mathbf{x}\\
      \mathbf{z}
      \end{bmatrix} = \mathcal{P}_{\{I, 0, 0, 0, G_1, G_2\}}\begin{bmatrix}
      \mathbf{x}\\
      \mathbf{z}_{ss}
      \end{bmatrix},
\end{equation}
\begin{equation}
\label{basic_identity_xs}
      \begin{bmatrix}
      \mathbf{x}\\
      \mathbf{z}_{s}
      \end{bmatrix} = \mathcal{P}_{\{I, 0, 0, 0, G_3, G_4\}}\begin{bmatrix}
      \mathbf{x}\\
      \mathbf{z}_{ss}
      \end{bmatrix},
\end{equation}
with
\begin{align}
G_1(s, \eta) &= B_a(s, \eta) + (s-\eta)I,\\
G_2(s, \eta) &= B_a(s, \eta),\\
G_3(s, \eta) &= B_b(s, \eta) + I,\\
G_4(s, \eta) &= B_b(s, \eta),\\
B_a(s, \eta) &= -\begin{bmatrix}
I & (s-a)I
\end{bmatrix}B_f(\eta),\\
B_b(\eta) &= -\begin{bmatrix}
0 & I
\end{bmatrix}B_f(\eta),\\
B_f(\eta) &= (B_c B_d)^{-1}B_c\begin{bmatrix}
0\\
(b-\eta)I\\
0\\
I
\end{bmatrix},\\
B_d &= \begin{bmatrix}
    I & 0\\
    I & (b-a)I\\
    0 & I\\
    0 & I
    \end{bmatrix}.
\end{align}
\end{lem}
\begin{pf}
For detailed proof, we refer to \citep[p. 2-3]{peet2018new}.
\QEDclosed
\end{pf}

\begin{lem}\label{total_fundamental_representation}
Suppose $\mathbf{x} \in \mathbb{R}^{n_{\text{o}}}$ and $\mathbf{z} \in W^{2,2}[a, b]$ and $$B_c\ \text{col}(\mathbf{z}(a), \mathbf{z}(b), \mathbf{z}_s(a), \mathbf{z}_s(b)) = 0,$$ with $B_c\in \R^{2n_{\text{p}}\times 4n_{\text{p}}}$ has full row rank. Then
\begin{equation}
    \label{combined-system_function}
    \begin{bmatrix}
\dot{\mathbf{x}}\\
\dot{\mathbf{z}}
    \end{bmatrix} = \begin{bmatrix}
    A+BL_1C & B L_2\mathcal{C}\\
    \mathcal{B}L_3C & \mathcal{A}+ \mathcal{B}L_4\mathcal{C}
    \end{bmatrix}\begin{bmatrix}
\mathbf{x}\\
\mathbf{z}
    \end{bmatrix}
\end{equation}
implies
\begin{equation}
    \label{combined-system_operator}
    \begin{bmatrix}
\dot{\mathbf{x}}\\
\dot{\mathbf{z}}
    \end{bmatrix} = \mathcal{P}_{\{U, V_1, V_2, Y, W_1, W_2\}}\begin{bmatrix}
\mathbf{x}\\
{\mathbf{z}}_{ss}
    \end{bmatrix},
\end{equation}
where
\begin{align*}
    &U = \Lambda_1 + \Lambda_2 , &Y = \Pi_1 + \Pi_2+A_2 + B_1 L_4 C_3, \nonumber\\
    &V_1 = \Omega_1 + \Omega_2 + BL_2C_3, \nonumber &W_1 = \Psi_1+\Psi_2, \nonumber\\
    &V_2 = \Sigma_1+ \Sigma_2,  &W_2 = \Upsilon_1 + \Upsilon_2, 
\end{align*}
with
\begin{align*}
& \big(\Lambda_1, \Omega_1, \Sigma_1, \Pi_1, \Psi_1, \Upsilon_1\big)\nonumber\\
&=\big(A+BL_1C, B L_2 C_a, B_1 L_3 C, A_0 + B_1 L_4 C_a, 0, 0\big)\nonumber\\ 
&\quad \qquad\times \big(I, 0, 0, 0, G_1, G_2\big), \\
    & \big(\Lambda_2, \Omega_2, \Sigma_2, \Pi_2, \Psi_2, \Upsilon_2\big)\nonumber \\ 
    &=\big(0, B L_2 C_b, 0, A_1 + B_1 L_4 C_b, 0, 0\big)\times \big(I, 0, 0, 0, G_3, G_4\big),
\end{align*}
\begin{align}
    C_3(s) = -C_1 B_d B_c(s) + C_1 \begin{bmatrix}
0\\
(b-s)I\\
0\\
I
\end{bmatrix}.
\end{align}
Moreover, the definitions of $G_1, G_2, G_3, G_4$ are given according to Lemma \ref{lemma_xs_x}.
\end{lem}
\begin{pf}
For detailed proof, see Appendix A.
\QEDclosed
\end{pf}
The above identities imply that the original state $\mathbf{z}$ can be fully reconstructed by the 'fundamental state' $\mathbf{z}_{ss}$ (see \cite{peet2018new}). Moreover, the new representation of the system dynamics in terms of $\mathbf{z}_{ss}$ is independent of the boundary constraints. In the subsequent sections, these identities are utilized for deriving stability conditions. 
\section{Lyapunov Stability Condition}
\label{section:stability}
In order to verify the stability of the PDE-ODE coupled system, we have proposed a Lyapunov function of the form
\begin{equation}
    \label{lyapunov_def}
    V(\mathbf{x},\mathbf{z}) = \Bigg\langle\begin{bmatrix}
\mathbf{x}\\
\mathbf{z}
    \end{bmatrix}, \mathcal{P}_{\{P,Q, Q^{\top}, S, R_1, R_2\}} \begin{bmatrix}
\mathbf{x}\\
\mathbf{z}
    \end{bmatrix} \Bigg\rangle_{X},
\end{equation}
such that $\mathcal{P}_{\{P,Q, Q^{\top}, S, R_1, R_2\}}$ is self-adjoint, i.e., $$\mathcal{P}_{\{P,Q, Q^{\top}, S, R_1, R_2\}} = \mathcal{P}_{\{P,Q, Q^{\top}, S, R_1, R_2\}}^{*}.$$
\subsection{Derivative of The Lyapunov Function}
\begin{thm}
    \label{LK_derivative}
Suppose $\mathbf{x} \in \mathbb{R}^{n_{\text{o}}}$ and $\mathbf{z} \in W^{2,2}[a, b]$ and $$B_c\ \text{col}(\mathbf{z}(a), \mathbf{z}(b), \mathbf{z}_s(a), \mathbf{z}_s(b)) = 0,$$ with $B_c\in \R^{2n_{\text{p}}\times 4n_{\text{p}}}$ has full row rank. Moreover, suppose there exists a matrix $P \in \mathbb{R}^{n_{\text{o}}\times n_{\text{o}}}$, and bounded polynomial functions $Q: [a, b] \rightarrow \mathbb{R}^{n_{\text{o}}\times n_{\text{p}}}$, $S: [a, b] \rightarrow \mathbb{R}^{n_{\text{p}}\times n_{\text{p}}}$, $R_1, R_2: [a, b]\times [a, b] \rightarrow \mathbb{R}^{n_{\text{p}}\times n_{\text{p}}}$, and
\begin{align}
    \bar{\mathbf{A}} = \begin{bmatrix}
    A+BL_1C & B L_2\mathcal{C}\\
    \mathcal{B}L_3C & \mathcal{A}+ \mathcal{B}L_4\mathcal{C}
    \end{bmatrix}.
\end{align}
Then 
\begin{align}
\label{derivative_expression}
      &\Bigg \langle  \begin{bmatrix}
    \mathbf{x}\\
    \mathbf{z}
    \end{bmatrix},{\mathcal{P}_{\{P,Q, Q^{\top}, S, R_1, R_2\}} \bar{\mbf{A}}\begin{bmatrix}
    \mathbf{x}\\
    \mathbf{z}
    \end{bmatrix}}\Bigg \rangle_{X}\\\nonumber
    &+ \Bigg \langle \bar{\mbf{A}}\begin{bmatrix}
    \mathbf{x}\\
    \mathbf{z}
    \end{bmatrix},{\mathcal{P}_{\{P,Q, Q^{\top}, S, R_1, R_2\}} \begin{bmatrix}
    \mathbf{x}\\
    \mathbf{z}
    \end{bmatrix}}\Bigg \rangle_{X}\\ 
    & = \Bigg\langle\begin{bmatrix}
\mathbf{x}\\
{\mathbf{z}}_{ss}
    \end{bmatrix}, \mathcal{P}_{\{\hat{K},\hat{L}, \hat{L}^{\top}, \hat{M}, \hat{N}_1, \hat{N}_2\}} \begin{bmatrix}
\mathbf{x}\\
{\mathbf{z}}_{ss}
    \end{bmatrix} \Bigg\rangle_{X}\nonumber
    \end{align}
where
\begin{align*}
    &\hat{K} = K + \bar{K}, &\hat{M} = M + \bar{M},  &\quad \hat{L} = L_1 + \bar{L}_1, \nonumber\\
    &\hat{L}^{\top} = L_2 + \bar{L}_2, &\hat{N}_1 = N_1 + \bar{N}_1,  &\quad \hat{N}_2 = N_2 + \bar{N}_2, \nonumber\\ 
\end{align*}
\begin{align}
    \big(\bar{K},\bar{L}_1, \bar{L}_2, \bar{M}, \bar{N}_1, \bar{N}_2\big) =& \big(K,L_1, L_2, M, N_1, N_2\big)^{*},\nonumber\\
    \big(K,L_1, L_2, M, N_1, N_2\big) =& \big(I,0, 0, 0, G_1, G_2\big)^{*} \nonumber\\
    &\times \big(P,Q, Q^{\top}, S, R_1, R_2\big) \nonumber\\
    &\times \big(U_1,V_1, V_2, Y, W_1, W_2\big).
\end{align}
Moreover, $G_1, G_2$ are as defined in Lemma \ref{lemma_xs_x} and $U$, $V_1$, $V_2$, $Y$, $W_1$, $W_2$ are as defined in Lemma \ref{total_fundamental_representation}. 

\end{thm}
\begin{pf}
The proof of Theorem \ref{LK_derivative} can be easily shown by substituting Lemma \ref{lemma_xs_x}-\ref{total_fundamental_representation} in \eqref{derivative_expression}.
\QEDclosed
\end{pf}
\subsection{Conditions of Exponential Stability}
In order to verify the exponential stability of the PDE-ODE coupled system, we search for the operator 
$$\mathcal{P}_{\{P, Q, Q^{\top}, S, R_1, R_2\}} \succ 0, \text{subject to}\ \mathcal{P}_{\{\hat{K}, \hat{L}, \hat{L}^{\top}, \hat{M}, \hat{N}_1, \hat{N}_2\}} \prec 0.$$
\begin{defn}
The cone of positive operators $\mathcal{P}$ that are given in the Definition \ref{def_operator} with polynomial multipliers and kernels associated with degree $d_1$ and $d_2$ are specified by
\begin{align}
    &\Phi_{d_1, d_2}:=\{(P, Q, Q^{\top}, S, R_1, R_2):\nonumber\\
&\text{the positivity condition is satisfied} \nonumber\\ 
& \quad \text{according to Theorem \ref{th:positivity}}\}.
\end{align}
\end{defn}
\begin{thm}
\label{thm:LMI}
   Suppose there exists $\epsilon >0, d_1, d_2 \in \mathbb{Z}$ with $P \in \mathbb{R}^{n_{\text{o}}}$, and bounded polynomial functions $Q: [a, b] \rightarrow \mathbb{R}^{n_{\text{o}}\times n_{\text{p}}}$, $S: [a, b] \rightarrow \mathbb{R}^{n_{\text{p}}\times n_{\text{p}}}$, $R_1, R_2: [a, b]\times [a, b] \rightarrow \mathbb{R}^{n_{\text{p}}\times n_{\text{p}}}$ such that
   \begin{align}
   \label{LMI_constraints}
     &\big(P-\epsilon I, Q, Q^{\top}, S-\epsilon I, R_1, R_2\big)\in \Phi_{d_1, d_2},\nonumber\\
     &-\big(\hat{K}, \hat{L}, \hat{L}^{\top}, M, N_1, N_2\big) - \epsilon \big(I, 0, 0, 0, T_1, T_2\big)\in \Phi_{d_1, d_2},
   \end{align}
   with
   \begin{align*}
    &\hat{K} = K + \bar{K}, &\hat{M} = M + \bar{M},  &\quad \hat{L} = L_1 + \bar{L}_1, \nonumber\\
    &\hat{L}^{\top} = L_2 + \bar{L}_2, &\hat{N}_1 = N_1 + \bar{N}_1,  &\quad \hat{N}_2 = N_2 + \bar{N}_2, \nonumber\\ 
\end{align*}
\begin{align}
    \big(\bar{K},\bar{L}_1, \bar{L}_2, \bar{M}, \bar{N}_1, \bar{N}_2\big) =& \big(K,L_1, L_2, M, N_1, N_2\big)^{*},\nonumber\\
    \big(K,L_1, L_2, M, N_1, N_2\big) =& \big(I,0, 0, 0, G_1, G_2\big)^{*} \nonumber\\
    &\times \big(P,Q, Q^{\top}, S, R_1, R_2\big) \nonumber\\
    &\times \big(U_1,V_1, V_2, Y, W_1, W_2\big), \nonumber\\
    \big(I,0, 0, 0, T_1, T_2\big) =& \big(I,0, 0, 0, G_1, G_2\big)^{*} \nonumber\\
    & \times \big(I,0, 0, 0, G_1, G_2\big),
\end{align}
where $G_1, G_2$ are as defined in Lemma \ref{lemma_xs_x} and $U$, $V_1$, $V_2$, $Y$, $W_1$, $W_2$ are as defined in Lemma \ref{total_fundamental_representation}.
Then the system in \eqref{combined-system} is exponentially stable. 

\end{thm}
\begin{pf}
We have defined the Lyapunov function as 
\begin{equation}
\begin{split}
    V(\mathbf{x},\mathbf{z}) =& \Bigg\langle\begin{bmatrix}
\mathbf{x}\\
\mathbf{z}
    \end{bmatrix}, \mathcal{P}_{\{P,Q, Q^{\top}, S, R_1, R_2\}} \begin{bmatrix}
\mathbf{x}\\
\mathbf{z}
    \end{bmatrix} \Bigg\rangle_{X}\\
   =&  \Bigg\langle\begin{bmatrix}
\mathbf{x}\\
\mathbf{z}
    \end{bmatrix}, \mathcal{P}_{\{P-\epsilon I,Q, Q^{\top}, S-\epsilon I, R_1, R_2\}} \begin{bmatrix}
\mathbf{x}\\
\mathbf{z}
    \end{bmatrix} \Bigg\rangle_{X}\\
    &+ \epsilon\Bigg\lvert\Bigg\lvert\begin{bmatrix}
\mathbf{x}\\
\mathbf{z}
    \end{bmatrix}\Bigg\rvert\Bigg\rvert_{X}^{2}\geq \epsilon\Bigg\lvert\Bigg\lvert\begin{bmatrix}
\mathbf{x}\\
\mathbf{z}
    \end{bmatrix}\Bigg\rvert\Bigg\rvert_{X}^{2}.
\end{split}
\end{equation}
This shows the strict positivity of \opPselfadj{P}{Q}{R}{S}. Furthermore, based on Theorem \ref{LK_derivative}, the exponential stability in Theorem \ref{th:stability_condition} can be reformulated as 
\begin{equation*}
    \begin{split}
    &\Bigg\langle\begin{bmatrix}
\mathbf{x}\\
{\mathbf{z}}_{ss}
    \end{bmatrix}, \mathcal{P}_{\{\hat{K},\hat{L}, \hat{L}^{\top}, \hat{N}_1, \hat{N}_2\}} \begin{bmatrix}
\mathbf{x}\\
{\mathbf{z}}_{ss}
    \end{bmatrix} \Bigg\rangle_{X}\\
    =& \Bigg\langle\begin{bmatrix}
\mathbf{x}\\
{\mathbf{z}}_{ss}
    \end{bmatrix}, \mathcal{P}_{\{\hat{K},\hat{L}, \hat{L}^{\top}, \hat{N}_1, \hat{N}_2\}} \begin{bmatrix}
\mathbf{x}\\
{\mathbf{z}}_{ss}
    \end{bmatrix} \Bigg\rangle_{X}\\
    &+ \epsilon\Bigg\lvert\Bigg\lvert\begin{bmatrix}
\mathbf{x}\\
\mathbf{z}
    \end{bmatrix}\Bigg\rvert\Bigg\rvert_{X}^{2} - \epsilon\Bigg\lvert\Bigg\lvert\begin{bmatrix}
\mathbf{x}\\
\mathbf{z}
    \end{bmatrix}\Bigg\rvert\Bigg\rvert_{X}^{2}\\
     =& \Bigg\langle\begin{bmatrix}
\mathbf{x}\\
{\mathbf{z}}_{ss}
    \end{bmatrix}, \mathcal{P}_{\{\hat{K},\hat{L}, \hat{L}^{\top}, \hat{N}_1, \hat{N}_2\}} \begin{bmatrix}
\mathbf{x}\\
{\mathbf{z}}_{ss}
    \end{bmatrix} \Bigg\rangle_{X}\\
    +& \epsilon \Bigg\langle
\mathcal{P}_{\{I,0, 0, 0, G_1, G_2\}}    \begin{bmatrix}
\mathbf{x}\\
{\mathbf{z}}_{ss}
    \end{bmatrix}, \mathcal{P}_{\{ I,0, 0, 0, G_1,  G_2\}} \begin{bmatrix}
\mathbf{x}\\
{\mathbf{z}}_{ss}
    \end{bmatrix} \Bigg\rangle_{X}\\
    &  - \epsilon\Bigg\lvert\Bigg\lvert\begin{bmatrix}
\mathbf{x}\\
\mathbf{z}
    \end{bmatrix}\Bigg\rvert\Bigg\rvert_{X}^{2},\\
    =& \Bigg\langle\begin{bmatrix}
\mathbf{x}\\
{\mathbf{z}}_{ss}
    \end{bmatrix}, \mathcal{P}_{\{\hat{K},\hat{L}_1, \hat{L}_2, \hat{M}, \hat{N}_1, \hat{N}_2\}} \begin{bmatrix}
\mathbf{x}\\
{\mathbf{z}}_{ss}
    \end{bmatrix} \Bigg\rangle_{X}\\
    +& \epsilon \Bigg\langle
\begin{bmatrix}
\mathbf{x}\\
{\mathbf{z}}_{ss}
    \end{bmatrix}, \mathcal{P}_{\{ I,0, 0, 0,  T_1,  T_2\}} \begin{bmatrix}
\mathbf{x}\\
{\mathbf{z}}_{ss}
    \end{bmatrix} \Bigg\rangle_{X}- \epsilon \Bigg\lvert\Bigg\lvert\begin{bmatrix}
\mathbf{x}\\
\mathbf{z}
    \end{bmatrix}\Bigg\rvert\Bigg\rvert_{X}^{2}\\
    &\leq - \epsilon \Bigg\lvert\Bigg\lvert\begin{bmatrix}
\mathbf{x}\\
\mathbf{z}
    \end{bmatrix}\Bigg\rvert\Bigg\rvert_{X}^{2}.
    \end{split}
\end{equation*}
This proves the exponential stability according to Theorem \ref{th:stability_condition}.
\QEDclosed
\end{pf}
In Theorem \ref{th:positivity}, we can choose $Z(s)$ to be a vector of monomials of degree $d_1$ and $Z(s, \eta)$ to be a vector of monomials of degree $d_2$. Hence, the constraint $(P, Q, Q^{\top}, S, R_1, R_2) \in \Phi_{d_1, d_2}$ can be viewed as an LMI constraint. As a result, the verification of stability for PDE-ODE coupled system amounts to a feasibility test of satisfying the LMI constraints related to \eqref{LMI_constraints} in Theorem \ref{thm:LMI}. 
\section{Numerical Implementations}
\label{section:results}
In this section, we illustrate the developed methodology for verifying the stability of PDE-ODE coupled system. The method has been implemented in MATLAB$^{\text{\textregistered}}$ using an adapted version of SOSTOOLS (see \cite{sostools}).
\subsection{ODE Coupled with Diffusion Equation}
First, we study the boundary controlled thermo-mechanical process where a lumped mechanical system is driven by converting thermal energy to mechanical work. In \cite{tang}, such a system is modeled as a finite dimensional ODE model with the actuator dynamics that are governed by the diffusion equation. A closed loop representation of such controlled PDE-ODE coupled system is given in \cite{tang} as
\begin{align}
    \dot{X}(t) &= -3X(t) + w(0, t),\\
w_t(x, t) &= w_{xx}(x t),\\
w_x(0, t) &= 0,\\
w(1, t) &= 0.
\end{align}
Based on the definitions in Section 3, we obtain $A = -3, B = 1, C=[]$ for the ODE system. For the PDE model, $A_0(s) = A_1(s) = [], A_2(s) = 1$. The inputs and output operators are $B_1(s) = [], C_1 = [1\ 0\ 0\ 0]$ with $C_a(s) = 0, C_b(s) = []$. The interconnection matrix $L = \begin{bmatrix}
0 & 1\\
1 & 0
\end{bmatrix}.$

The developed methodology proves the exponential stability of this system which has been verified analytically in \cite{tang}.

\subsection{ODE Coupled with Diffusion-Reaction PDE}
The next example is taken from chemical processes where due to a clear distinction between microscopic processes and macroscopic processes, models often involve PDEs governing macroscopic dynamics coupled with ODEs governing microscopic dynamics. For exothermic or endothermic reaction processes, we consider a general class of forced diffusion-reaction in the spatial domain $[0,1]$ which is represented by the following coupled PDEs.
\begin{equation}
    \label{Ahmadi}
    \dot{z}(s,t) = \lambda z(s,t) + z_{ss}(s,t) + B_1(s) u(t).
    \end{equation}
In \cite{ahmadi1}, it has been shown that if the boundary conditions are chosen to be $z(0,t) = z(1,t) = 0$, then the unforced system given in \eqref{Ahmadi}(i.e. $B_1(s) = 0$) is stable if and only if $\lambda < \pi^2 = 9.8696$. Similarly, in \cite{valmorbida2016stability} it has been shown that for the boundary conditions $z(0,t) = z_s(1,t) = 0$, the system is stable if and only if $\lambda \leq 2.467$. 

Such diffusion-reaction processes in \eqref{Ahmadi} are coupled with a finite dimensional system that typically represents a molecular transport in the microscopic scale. Specifically, we consider that the effect of these molecular transport applies a uniform in-domain input to the PDEs over the entire spatial domain (i.e. $B_1(s) = 1, \forall s \in [a, b]$). On the other hand, the output flux of the chemical reaction at $s =0$(i.e. $z_{s}(0,t) = 0$) drives the molecular transport and is considered to be the input to the ODE model. The interconnection matrix $L = \begin{bmatrix}
0 & I\\
I & 0
\end{bmatrix}$. The ODE model is given in \eqref{ode_ahmadi}. The stability of \eqref{ode_ahmadi} can be verified by checking the eigenvalues. 
\begin{equation}
\label{ode_ahmadi}
    \begin{split}
       \dot{x}(t) &= \begin{bmatrix}
       -1.2142  &  1.9649  &  0.2232 &   0.5616\\
   -1.8042 &  -0.7260 &  -0.3479  &  5.4355\\
   -0.2898  &  0.7381  & -1.7606  &  0.8294\\
   -0.9417  & -5.3399  & -1.0704  & -0.7590
       \end{bmatrix}{x}(t)\\
       &+ \begin{bmatrix}
         -1.5368    &     0\\
         0 &   0.8871\\
    1.0656   &      0\\
    1.1882    &     0
       \end{bmatrix}v(t),\\
       w(t) &= \begin{bmatrix}
       -2.5575    &     0  &  1.0368     &    0\\
   -1.8067  &  0.4630  &  1.3621     &    0
       \end{bmatrix}{x}(t).
    \end{split}
\end{equation}
The developed methodology is able to verify the stability of the coupled system for a) $\lambda < \pi^2 = 9.8696$ with boundary conditions $z(0,t) = z(1,t) = 0$, and b) $\lambda \leq 2.467$ with boundary conditions $z(0,t) = z_s(1,t) = 0$.

\section{Conclusions}
\label{section:conclusions}
In this paper, we have presented a novel methodology to verify stability for a large class of mutually coupled linear PDE-ODE systems in one spatial dimension. In particular, by using positive-definite matrices to parametrize a cone of Lyapunov functions, we have shown that the stability of linear PDE-ODE coupled systems amounts to satisfying a specific Linear Matrix Inequality (LMI) that can be efficiently solved using convex optimization. We have provided a scalable variation of Sum of Squares (SOS) algorithm that can be utilized to verify the stability a large class of PDE-ODE coupled system efficiently.

In future, the results of stability analysis for PDE-ODE coupled system can be extended to study the robustness properties and $\mathcal{H}_{\infty}$ analysis of PDE-ODE coupled system. Extending the current framework for PDEs in more than one spatial dimension and nonlinear PDEs are a potential direction of future research. Furthermore, viewing the coupling of ODEs and PDEs as a feedback interconnection, the problem of synthesizing finite dimensional controller and observer are of great practical importance.
\bibliography{ifacconf.bbl} 

\appendix
\section{Proof of Lemma \ref{total_fundamental_representation}}
\textit{Lemma \ref{total_fundamental_representation}. }
Suppose $\mathbf{x} \in \mathbb{R}^{n_{\text{o}}}$ and $\mathbf{z} \in W^{2,2}[a, b]$ and $$B_c\ \text{col}(\mathbf{z}(a), \mathbf{z}(b), \mathbf{z}_s(a), \mathbf{z}_s(b)) = 0,$$ with $B_c\in \R^{2n_{\text{p}}\times 4n_{\text{p}}}$ has full row rank. Then
\begin{equation}
    \label{combined-system_function_pf}
    \begin{bmatrix}
\dot{\mathbf{x}}\\
\dot{\mathbf{z}}
    \end{bmatrix} = \begin{bmatrix}
    A+BL_1C & B L_2\mathcal{C}\\
    \mathcal{B}L_3C & \mathcal{A}+ \mathcal{B}L_4\mathcal{C}
    \end{bmatrix}\begin{bmatrix}
\mathbf{x}\\
\mathbf{z}
    \end{bmatrix}
\end{equation}
implies
\begin{equation}
    \label{combined-system_operator_pf}
    \begin{bmatrix}
\dot{\mathbf{x}}\\
\dot{\mathbf{z}}
    \end{bmatrix} = \mathcal{P}_{\{U, V_1, V_2, Y, W_1, W_2\}}\begin{bmatrix}
\mathbf{x}\\
{\mathbf{z}}_{ss}
    \end{bmatrix},
\end{equation}
where
\begin{align*}
    &U = \Lambda_1 + \Lambda_2 , &Y = \Pi_1 + \Pi_2+A_2 + B_1 L_4 C_3, \nonumber\\
    &V_1 = \Omega_1 + \Omega_2 + BL_2C_3, \nonumber &W_1 = \Psi_1+\Psi_2, \nonumber\\
    &V_2 = \Sigma_1+ \Sigma_2,  &W_2 = \Upsilon_1 + \Upsilon_2, 
\end{align*}
with
\begin{align*}
& \big(\Lambda_1, \Omega_1, \Sigma_1, \Pi_1, \Psi_1, \Upsilon_1\big)\nonumber\\
&=\big(A+BL_1C, B L_2 C_a, B_1 L_3 C, A_0 + B_1 L_4 C_a, 0, 0\big)\nonumber\\ 
&\quad \qquad\times \big(I, 0, 0, 0, G_1, G_2\big), \\
    & \big(\Lambda_2, \Omega_2, \Sigma_2, \Pi_2, \Psi_2, \Upsilon_2\big)\nonumber \\ 
    &=\big(0, B L_2 C_b, 0, A_1 + B_1 L_4 C_b, 0, 0\big)\times \big(I, 0, 0, 0, G_3, G_4\big),
\end{align*}
\begin{align}
    C_3(s) = -C_1 B_d B_c(s) + C_1 \begin{bmatrix}
0\\
(b-s)I\\
0\\
I
\end{bmatrix}.
\end{align}
Moreover, the definitions of $G_1, G_2, G_3, G_4$ are given according to Lemma \ref{lemma_xs_x}.
\begin{pf}
First of all, using \eqref{PDE_boundary} we can show that
\begin{align}
C_1 \mathbf{z}_b = \int\limits_{a}^b C_3(s)\mathbf{z}_{ss}(s) \text{d}s,
\end{align}
where 
\begin{align}
    C_3(s) = -C_1 B_d B_c(s) + C_1 \begin{bmatrix}
0\\
(b-s)I\\
0\\
I
\end{bmatrix}.
\end{align}
Now we use the definitions $\mathcal{A}, \mathcal{B}, \mathcal{C}$ in \eqref{PDE_dynamics} and re-write \eqref{combined-system_function} as
\begin{align}
  \begin{bmatrix}
\dot{\mathbf{x}}\\
\dot{\mathbf{z}}
    \end{bmatrix} =& \begin{bmatrix}
    A+BL_1C & BL_2\mathcal{C}_a\\
    \mathcal{B}L_3C & \mathcal{A}_0+ \mathcal{B}L_4\mathcal{C}_a
    \end{bmatrix}\begin{bmatrix}
\mathbf{x}\\
\mathbf{z}
    \end{bmatrix}\nonumber\\
    &+  \begin{bmatrix}
    0 & B L_2 \mathcal{C}_b\\
    0 & \mathcal{A}_1+ \mathcal{B}L_4\mathcal{C}_b
    \end{bmatrix}\begin{bmatrix}
\mathbf{x}\\
\mathbf{z}_{s}
    \end{bmatrix} + \begin{bmatrix}
    0 & B L_2 \mathcal{C}_3\\
    0 & \mathcal{A}_2+ \mathcal{B}L_4\mathcal{C}_3
    \end{bmatrix}\begin{bmatrix}
\mathbf{x}\\
\mathbf{z}_{ss}
    \end{bmatrix}.
\end{align}
Based on \eqref{PDE_dynamics},

$[\mathcal{A}_0\mathbf{z}](s) := A_0(s)\mathbf{z}(s), [\mathcal{A}_1\mathbf{z}_s](s) := A_1(s)\mathbf{z}_s(s)$,

$[\mathcal{A}_2\mathbf{z}_{ss}](s) := A_2(s)\mathbf{z}_{ss}(s), [\mathcal{C}_a\mathbf{z}](s):=\int\limits_{a}^{b} C_a(s)\mathbf{z}(s)\text{d}s$,

$[\mathcal{C}_b\mathbf{z}_{s}](s):=\int\limits_{a}^{b} C_b(s)\mathbf{z}_{s}(s)\text{d}s$, 

and $[\mathcal{C}_3\mathbf{z}_{ss}](s):=\int\limits_{a}^{b} C_3(s)\mathbf{z}_{ss}(s)\text{d}s$.

Rewriting in the form of \opPgen
\begin{align}
   \begin{bmatrix}
\dot{\mathbf{x}}\\
\dot{\mathbf{z}}
    \end{bmatrix} =& \mathcal{P}_{\{A+BL_1C, B L_2 C_a, B_1 L_3 C, A_0 + B_1 L_4 C_a, 0, 0\}}\begin{bmatrix}
\mathbf{x}\\
\mathbf{z}
    \end{bmatrix}\nonumber\\
    &+ \mathcal{P}_{\{0, B L_2 C_b, 0, A_1 + B_1 L_4 C_b, 0, 0\}} \begin{bmatrix}
\mathbf{x}\\
\mathbf{z}_{s}
    \end{bmatrix}\nonumber\\
    &+ \mathcal{P}_{\{0, B L_2 C_3, 0, A_2 + B_1 L_4 C_3, 0, 0\}}\begin{bmatrix}
\mathbf{x}\\
\mathbf{z}_{ss}
    \end{bmatrix}.   
\end{align}
Now, substituting Lemma \ref{lemma_xs_x} in \eqref{lemma_system:eq2} we obtain
\begin{align}
\label{lemma_system:eq2}
   \begin{bmatrix}
\dot{\mathbf{x}}\\
\dot{\mathbf{z}}
    \end{bmatrix} =& \mathcal{P}_{\{\Lambda_1, \Omega_1, \Sigma_1, \Pi_1, \Psi_1, \Upsilon_1\}}\begin{bmatrix}
\mathbf{x}\\
\mathbf{z}_{ss}
    \end{bmatrix}\nonumber\\
    &+ \mathcal{P}_{\{\Lambda_2, \Omega_2, \Sigma_2, \Pi_2, \Psi_2, \Upsilon_2\}} \begin{bmatrix}
\mathbf{x}\\
\mathbf{z}_{ss}
    \end{bmatrix}\nonumber\\
    &+ \mathcal{P}_{\{0, B L_2 C_3, 0, A_2 + B_1 L_4 C_3, 0, 0\}}\begin{bmatrix}
\mathbf{x}\\
\mathbf{z}_{ss}
    \end{bmatrix},   
\end{align}
where $\mathcal{P}_{\{\Lambda_1, \Omega_1, \Sigma_1, \Pi_1, \Psi_1, \Upsilon_1\}}$ is defined as
\begin{align*}
\mathcal{P}_{\{A+BL_1C, B L_2 C_a, B_1 L_3 C, A_0 + B_1 L_4 C_a, 0, 0\}} \mathcal{P}_{\{I, 0, 0, 0, G_1, G_2\}}.
\end{align*}
Similarly, $\mathcal{P}_{\{\Lambda_2, \Omega_2, \Sigma_2, \Pi_2, \Psi_2, \Upsilon_2\}}$ is defined as
\begin{align*}
\mathcal{P}_{\{0, B L_2 C_b, 0, A_1 + B_1 L_4 C_b, 0, 0\}} \mathcal{P}_{\{I, 0, 0, 0, G_3, G_4\}}.
\end{align*}
As a result, 
\begin{align*}
&\mathcal{P}_{\{\Lambda_1, \Omega_1, \Sigma_1, \Pi_1, \Psi_1, \Upsilon_1\}} +  \mathcal{P}_{\{\Lambda_2, \Omega_2, \Sigma_2, \Pi_2, \Psi_2, \Upsilon_2\}}\nonumber\\
&+\mathcal{P}_{\{0, B L_2 C_3, 0, A_2 + B_1 L_4 C_3, 0, 0\}} = \mathcal{P}_{\{U, V_1, V_2, Y, W_1, W_2\}}.\QEDclosed
\end{align*}
\end{pf}

\section{Proof of Theorem \ref{LK_derivative}}
\textit{Theorem \ref{LK_derivative}. }
Suppose $\mathbf{x} \in \mathbb{R}^{n_{\text{o}}}$ and $\mathbf{z} \in W^{2,2}[a, b]$ and $$B_c\ \text{col}(\mathbf{z}(a), \mathbf{z}(b), \mathbf{z}_s(a), \mathbf{z}_s(b)) = 0,$$ with $B_c\in \R^{2n_{\text{p}}\times 4n_{\text{p}}}$ has full row rank. Moreover, suppose there exists a matrix $P \in \mathbb{R}^{n_{\text{o}}\times n_{\text{o}}}$, and bounded polynomial functions $Q: [a, b] \rightarrow \mathbb{R}^{n_{\text{o}}\times n_{\text{p}}}$, $S: [a, b] \rightarrow \mathbb{R}^{n_{\text{p}}\times n_{\text{p}}}$, $R_1, R_2: [a, b]\times [a, b] \rightarrow \mathbb{R}^{n_{\text{p}}\times n_{\text{p}}}$, and
\begin{align}
    \bar{\mathbf{A}} = \begin{bmatrix}
    A+BL_1C & B L_2\mathcal{C}\\
    \mathcal{B}L_3C & \mathcal{A}+ \mathcal{B}L_4\mathcal{C}
    \end{bmatrix}.
\end{align}
Then 
\begin{align}
\label{derivative_expression_pf}
      &\Bigg \langle  \begin{bmatrix}
    \mathbf{x}\\
    \mathbf{z}
    \end{bmatrix},{\mathcal{P}_{\{P,Q, Q^{\top}, S, R_1, R_2\}} \bar{\mbf{A}}\begin{bmatrix}
    \mathbf{x}\\
    \mathbf{z}
    \end{bmatrix}}\Bigg \rangle_{X}\\\nonumber
    &+ \Bigg \langle \bar{\mbf{A}}\begin{bmatrix}
    \mathbf{x}\\
    \mathbf{z}
    \end{bmatrix},{\mathcal{P}_{\{P,Q, Q^{\top}, S, R_1, R_2\}} \begin{bmatrix}
    \mathbf{x}\\
    \mathbf{z}
    \end{bmatrix}}\Bigg \rangle_{X}\\ 
    & = \Bigg\langle\begin{bmatrix}
\mathbf{x}\\
{\mathbf{z}}_{ss}
    \end{bmatrix}, \mathcal{P}_{\{\hat{K},\hat{L}, \hat{L}^{\top}, \hat{M}, \hat{N}_1, \hat{N}_2\}} \begin{bmatrix}
\mathbf{x}\\
{\mathbf{z}}_{ss}
    \end{bmatrix} \Bigg\rangle_{X}\nonumber
    \end{align}
where
\begin{align*}
    &\hat{K} = K + \bar{K}, &\hat{M} = M + \bar{M},  &\quad \hat{L} = L_1 + \bar{L}_1, \nonumber\\
    &\hat{L}^{\top} = L_2 + \bar{L}_2, &\hat{N}_1 = N_1 + \bar{N}_1,  &\quad \hat{N}_2 = N_2 + \bar{N}_2, \nonumber\\ 
\end{align*}
\begin{align}
    \big(\bar{K},\bar{L}_1, \bar{L}_2, \bar{M}, \bar{N}_1, \bar{N}_2\big) =& \big(K,L_1, L_2, M, N_1, N_2\big)^{*},\nonumber\\
    \big(K,L_1, L_2, M, N_1, N_2\big) =& \big(I,0, 0, 0, G_1, G_2\big)^{*} \nonumber\\
    &\times \big(P,Q, Q^{\top}, S, R_1, R_2\big) \nonumber\\
    &\times \big(U_1,V_1, V_2, Y, W_1, W_2\big).
\end{align}
Moreover, $G_1, G_2$ are as defined in Lemma \ref{lemma_xs_x} and $U$, $V_1$, $V_2$, $Y$, $W_1$, $W_2$ are as defined in Lemma \ref{total_fundamental_representation}.
\begin{pf}
The time derivative of the Lyapunov function is
\begin{equation}
\begin{split}
    &\Bigg \langle  \begin{bmatrix}
    \mathbf{x}\\
    \mathbf{z}
    \end{bmatrix},{\mathcal{P}_{\{P,Q, Q^{\top}, S, R_1, R_2\}} \bar{\mbf{A}}\begin{bmatrix}
    \mathbf{x}\\
    \mathbf{z}
    \end{bmatrix}}\Bigg \rangle_{X}\\
    &+ \Bigg \langle \bar{\mbf{A}}\begin{bmatrix}
    \mathbf{x}\\
    \mathbf{z}
    \end{bmatrix},{\mathcal{P}_{\{P,Q, Q^{\top}, S, R_1, R_2\}} \begin{bmatrix}
    \mathbf{x}\\
    \mathbf{z}
    \end{bmatrix}}\Bigg \rangle_{X}\\ 
    &= \underbrace{\Bigg\langle\mathcal{P}_{\{U_1,V_1, V_2, Y, W_1, W_2\}} \begin{bmatrix}
\mathbf{x}\\
{\mathbf{z}}_{ss}
    \end{bmatrix}, \mathcal{P}_{\{P,Q, Q^{\top}, S, R_1, R_2\}} \begin{bmatrix}
\mathbf{x}\\
\mathbf{z}
    \end{bmatrix} \Bigg\rangle_{X}}_{:= \Gamma_1}\\
    &+\underbrace{\Bigg\langle\begin{bmatrix}
\mathbf{x}\\
\mathbf{z}
    \end{bmatrix}, \mathcal{P}_{\{P,Q, Q^{\top}, S, R_1, R_2\}} \mathcal{P}_{\{U_1,V_1, V_2, Y, W_1, W_2\}}\begin{bmatrix}
\mathbf{x}\\
{\mathbf{z}}_{ss}
    \end{bmatrix} \Bigg\rangle_{X}}_{:= \Gamma_2}.
    \end{split}
\end{equation}
Now, taking two individual terms separately and substituting the identity \eqref{basic_identity_x} we obtain
\begin{equation}
\begin{split}
    \Gamma_2&= \Bigg\langle \begin{bmatrix}
\mathbf{x}\\
{\mathbf{z}}_{ss}
    \end{bmatrix},\mathcal{P}_{\{I,0, 0, 0, G_1, G_2\}}^{*}\\
    &\mathcal{P}_{\{P,Q, Q^{\top}, S, R_1, R_2\}} \mathcal{P}_{\{U_1,V_1, V_2, Y, W_1, W_2\}} \begin{bmatrix}
\mathbf{x}\\
{\mathbf{z}}_{ss}
    \end{bmatrix} \Bigg\rangle_{X}.\\
    & = \Bigg\langle\begin{bmatrix}
\mathbf{x}\\
{\mathbf{z}}_{ss}
    \end{bmatrix}, \mathcal{P}_{\{K,L_1, L_2,M, N_1, N_2\}} \begin{bmatrix}
\mathbf{x}\\
{\mathbf{z}}_{ss}
    \end{bmatrix} \Bigg\rangle_{X}.
\end{split}
\end{equation}
Similarly
\begin{equation}
\begin{split}
    \Gamma_1&= \Bigg\langle\mathcal{P}_{\{U_1,V_1, V_2, Y, W_1, W_2\}} \begin{bmatrix}
\mathbf{x}\\
{\mathbf{z}}_{ss}
    \end{bmatrix},\\
    &\mathcal{P}_{\{P,Q, Q^{\top}, S, R_1, R_2\}} \mathcal{P}_{\{I,0, 0, 0, G_1, G_2\}} \begin{bmatrix}
\mathbf{x}\\
{\mathbf{z}}_{ss}
    \end{bmatrix} \Bigg\rangle_{X}.\\
   &= \Bigg\langle \begin{bmatrix}
\mathbf{x}\\
{\mathbf{z}}_{ss}
    \end{bmatrix},\mathcal{P}_{\{U_1,V_1, V_2, Y, W_1, W_2\}}^{*}\\
    &\mathcal{P}_{\{P,Q, Q^{\top}, S, R_1, R_2\}} \mathcal{P}_{\{I,0, 0, 0, G_1, G_2\}} \begin{bmatrix}
\mathbf{x}\\
{\mathbf{z}}_{ss}
    \end{bmatrix} \Bigg\rangle_{X}.\\
    & = \Bigg\langle\begin{bmatrix}
\mathbf{x}\\
{\mathbf{z}}_{ss}
    \end{bmatrix}, \mathcal{P}_{\{K,L_1, L_2,M, N_1, N_2\}}^{*} \begin{bmatrix}
\mathbf{x}\\
{\mathbf{z}}_{ss}
    \end{bmatrix} \Bigg\rangle_{X}
\end{split}
\end{equation}
Now, using the adjoint operation of the operators, we obtain
\begin{equation}
    \begin{split}
    &\Bigg\langle\begin{bmatrix}
\mathbf{x}\\
{\mathbf{z}}_{ss}
    \end{bmatrix}, \mathcal{P}_{\{K,L_1, L_2, M, N_1, N_2\}} \begin{bmatrix}
\mathbf{x}\\
{\mathbf{z}}_{ss}
    \end{bmatrix} \Bigg\rangle_{X}\\
    &+\Bigg\langle\begin{bmatrix}
\mathbf{x}\\
{\mathbf{z}}_{ss}
    \end{bmatrix}, \mathcal{P}_{\{K,L_1, L_2, M, N_1, N_2\}}^{*} \begin{bmatrix}
\mathbf{x}\\
{\mathbf{z}}_{ss}
    \end{bmatrix} \Bigg\rangle_{X}\\    
    &= \Bigg\langle\begin{bmatrix}
\mathbf{x}\\
{\mathbf{z}}_{ss}
    \end{bmatrix},\mathcal{P}_{\{\hat{K},\hat{L}, \hat{L}^{\top},\hat{M}, \hat{N}_1, \hat{N}_2\}}\begin{bmatrix}
\mathbf{x}\\
{\mathbf{z}}_{ss}
    \end{bmatrix} \Bigg\rangle_{X}.
    \end{split}
\end{equation}
with 
\begin{align*}
\mathcal{P}_{\{\hat{K},\hat{L}, \hat{L}^{\top},\hat{M}, \hat{N}_1, \hat{N}_2\}} =& \mathcal{P}_{\{K,L_1, L_2,M, N_1, N_2\}}\\ 
&+\qquad \mathcal{P}_{\{K,L_1, L_2,M, N_1, N_2\}}^{*},
\end{align*}
This completes the proof. 
\QEDclosed
\end{pf}
\section{Composition of the Operators}
\textit{Lemma \ref{lem:composition}. } For any matrices $A,P\in\R^{m\times m}$ and bounded functions $B_1, Q_1: [a, b]\to\R^{m\times n}$, $B_2, Q_2: [a, b]\to\R^{n\times m}$, $D,S:[a, b]\to\R^{n \times n}$, $C_i,R_i:[a, b]\times [a, b] \to \R^{n\times n}$ with $i \in \{1,2\}$, the following identity holds 
\begin{align*}
    \opPpq{A}{B}{C}{D}&\opPpq{P}{Q}{R}{S}\\
    &\hspace{-1cm}= \mc{P}_{\{\hat{P},\hat{Q}_1,\hat{Q}_2,\hat{S},\hat{R}_1,\hat{R}_2\}}
\end{align*}
where
\begin{align*}
    \hat{P} &= AP + \myint B_1(s)Q_2(s)\text{d}s,\\
    \hat{Q}_1(s) &= AQ_1(s) + B_1(s)S(s)+\myintb{s}B_1(\eta)R_1(\eta,s)\text{d}\eta\\
    &\qquad+\myinta{s}B_1(\eta)R_2(\eta,s)\text{d}\eta,\\
    \hat{Q}_2(s) &= B_2(s)P + D(s)Q_2(s) + \myinta{s}C_1(s,\eta)Q_2(\eta)\text{d}\eta\\
    &\qquad+\myintb{s}C_2(s,\eta)Q_2(\eta)\text{d}\eta,\\
    \hat{S}(s) &= D(s)S(s),\\
    \hat{R}_1(s,\eta) &=B_2(s)Q_1(\eta)+D(s)R_1(s,\eta)+C_1(s,\eta)S(\eta)\\
    &\hspace{-0.5cm}+\myinta{\eta} C_1(s,\theta)R_2(\theta,\eta)\text{d}\theta+\int_{\eta}^{s}C_1(s,\theta)R_1(\theta,\eta)\text{d}\theta\\
    &\hspace{-0.5cm}+\myintb{s}C_2(s,\theta)R_1(\theta,\eta)\text{d}\theta,\\
    \hat{R}_2(s,\eta) &=B_2(s)Q_1(\eta)+D(s)R_2(s,\eta)+C_2(s,\eta)S(\eta)\\
    &\hspace{-0.5cm}+\myinta{s} C_1(s,\theta)R_2(\theta,\eta)\text{d}\theta+\int_{s}^{\eta}C_2(s,\theta)R_2(\theta,\eta)d\theta\\
    &\hspace{-0.5cm}+\myintb{\eta}C_2(s,\theta)R_1(\theta,\eta)\text{d}\theta.
\end{align*}
\begin{pf}
Suppose
\begin{align*}
    &\opPpq{A}{B}{C}{D}\left(\opPpq{P}{Q}{R}{S}\vmatwo{x_1}{x_2}\right)(s) \\
    &\qquad= \left(\mc{P}_{\{\hat{P},\hat{Q}_1,\hat{Q}_2,\hat{S},\hat{R}_1,\hat{R}_2\}}\vmatwo{x_1}{x_2}\right)(s)
\end{align*}
Let
\begin{align*}
    \left(\opPpq{P}{Q}{R}{S}\vmatwo{x_1}{x_2}\right)(s)= \vmatwo{y_1}{y_2(s)}
\end{align*}
where $y_1 = Px_1 + \myint Q_1(s) x_2(s)\text{d}s$ and $y_2(s)= Q_2(s)x_1 + S(s)x_2(s)+\myinta{s}R_1(s,\eta)x_2(\eta)\text{d}\eta+\myintb{s}R_2(s,\eta)x_2(\eta)\text{d}\eta$.

Then,
\begin{align*}
    \opPpq{A}{B}{C}{D}\opPpq{P}{Q}{R}{S}&\vmatwo{x_1}{x_2}&\\
    &\hspace{-4cm}=\left(\opPpq{A}{B}{C}{D}\vmatwo{y_1}{y_2}\right)(s)= \vmatwo{z_1}{z_2(s)}
\end{align*}
where $z_1 = Ay_1 + \myint B_1(s) y_2(s)\text{d}s$ and $z_2(s)= B_2(s)y_1 + D(s)y_2(s)+\myinta{s}C_1(s,\eta)y_2(\eta)\text{d}\eta+\myintb{s}C_2(s,\eta)y_2(\eta)\text{d}\eta$.

Finding the composition is a straight-forward algebraic operation. We do this by expanding each term separately. Firstly,
\begin{align*}
    &\myint B_1(s) y_2(s) \text{d}s = \myint B_1(s) \Big(Q_2(s)x_1 + S(s)x_2(s)\\
    &~~+\myinta{s}R_1(s,\eta)x_2(\eta)\text{d}\eta+\myintb{s}R_2(s,\eta)x_2(\eta)\text{d}\eta\Big)\text{d}s.\\
    &\hspace{-0.35cm}\text{Seperating the terms involving $x_1$ and $x_2(s)$, we obtain}\\
    &z_1 = \underbrace{\left(AP+\myint B_1(s)Q_2(s)\right)}_{\hat{P}}x_1+\myint \hat{Q}_1(s) x_2(s) ds
\end{align*}
where
\begin{align*}
    \hat{Q}_1(s) &= AQ_1(s) + B_1(s)S(s)+\myintb{s}B_1(\eta)R_1(\eta,s)\text{d}\eta\\
    &\qquad+\myinta{s}B_1(\eta)R_2(\eta,s)\text{d}\eta
\end{align*}

Next, we expand the terms of $z_2(s)$ and specifically collect the terms having $x_1$. We obtain
\begin{align*}
    \Big(B_2(s)P + D(s)Q_2(s) &+ \myinta{s}C_1(s,\eta)Q_2(\eta)\text{d}\eta\Big)x_1 \\
    &= \hat{Q}_2(s) x_1
\end{align*}
Next, grouping the terms having $x_2(s)$ we obtain
\begin{align*}
    D(s)S(s) x_2(s) = \hat{S}(s) x_2(s).
\end{align*}

Similarly, we can group the terms involving $x_2(\th)$ as
\begin{align*}
    \myinta{s} &\Bigg(B_2(s)Q_1(\eta)+D(s)R_1(s,\eta)+C_1(s,\eta)S(\eta)\\
    &+\int_{\eta}^{s}C_1(s,\theta)R_1(\theta,\eta)\text{d}\theta+\myinta{\eta} C_1(s,\theta)R_2(\theta,\eta)\text{d}\theta\\
    &+\myintb{s}C_2(s,\theta)R_1(\theta,\eta)d\theta\Bigg) x_2(\eta) \text{d}\eta \\
    &=\myinta{s} \hat{R_1}(s,\eta)x_2(\eta)\text{d}\eta. 
\end{align*}    
and 
\begin{align*}
     \myintb{s}&\Bigg(B_2(s)Q_1(\eta)+D(s)R_2(s,\eta)+C_2(s,\eta)S(\eta)\\
    &+\int_{s}^{\eta}C_2(s,\theta)R_2(\theta,\eta)\text{d}\theta+\myinta{s} C_1(s,\theta)R_2(\theta,\eta)\text{d}\theta\\
    &+\myintb{\th}C_2(s,\theta)R_1(\theta,\eta)\text{d}\theta\Bigg)x_2(\eta)\text{d}\eta \\
    &= \myintb{s} \hat{R_2}(s,\eta)x_2(\eta) \text{d}\eta.
\end{align*}
This completes the proof.
\QEDclosed
\end{pf}
\section{Adjoint of the Operators}
\textit{Lemma \ref{lem:adjoint}. } For any matrices $P\in\R^{m\times m}$ and  bounded functions $Q_1: [a, b]\to\R^{m\times n}$, $Q_2: [a, b]\to\R^{n\times m}$, $S:[a, b]\to\R^{n \times n}$, $R_1, R_2:[a, b]\times [a, b] \to \R^{n\times n}$,  the following identity holds for any $x \in \mathbb{R}^{m}, y \in L_2^n([a, b])$.
\begin{align}
    &\Big\langle{x},{\opPpq{P}{Q}{R}{S}y}\Big\rangle_{\R^{m}\times L_2^n[a, b]}\nonumber\\
    &=\Big\langle{\opPpq{P}{Q}{R}{S}^*x},{y}\Big\rangle_{\R^{m}\times L_2^n[a, b]},
\end{align}
where, $\opPpq{P}{Q}{R}{S}^*=\opPpq{\hat{P}}{\hat{Q}}{\hat{R}}{\hat{S}}$ with
\begin{align}
\label{adjoint_matrix_pf}
    &\hat{P} = P^{\top}, &\hat{S}(s) = S^{\top}(s), \nonumber\\
    &\hat{Q}_1(s) = Q_2^{\top}(s), \nonumber &\hat{R}_1(s,\eta) = R_2^{\top}(\eta,s), \nonumber\\
    &\hat{Q}_2(s) = Q_1^{\top}(s),  &\hat{R}_2(s,\eta) = R_1^{\top}(\eta,s). 
\end{align}

\begin{pf}
We use the fact that for any scalar $a$ we have $a=a^{\top}$. Let  $ x=\vmatwo{x_1}{x_2(s)}$ and $y=\vmatwo{y_1}{y_2(s)}$. 


Then {\smaller
\begin{align*}
    &\Big\langle{x},{\opPpq{P}{Q}{R}{S}y}\Big\rangle_{\R^{m}\times L_2^n[a, b]}\\
    &~~= x_1^{\top}Py_1+\myint x_1^{\top}Q_1(s)y_2(s)\text{d}s+\myint y_2^{\top}(s) Q_2(s) x_1 \text{d}s \\
    &~~~+ \myint x_2(s)^{\top} S(s) y_2(s) \text{d}s + \myint\myinta{s}x_2^{\top}(s) R_1(s,\eta)y_2(\eta) \text{d}\eta \text{d}s\\
    &~~~+\myint\myintb{s}x_2^{\top}(s) R_2(s,\eta)y_2(\eta) \text{d}\eta \text{d}s\\
    &~~= y_1^{\top}P^{\top}x_1+\myint x_1^{\top}Q_2^{\top}(s)y_2(s)\text{d}s+\myint y_2(s) Q_1^{\top}(s) x_1 \text{d}s \\
    &~~~+ \myint y_2^{\top}(s) S^{\top}(s) x_2(s) \text{d}s + \myint\myinta{s}y_2^{\top}(s) R_2^{\top}(\eta,s)x_2(\eta) \text{d}\eta \text{d}s\\
    &~~~+\myint\myintb{s}y_2^{\top}(s) R_1^{\top}(\eta,s)x_2(\eta) \text{d}\eta \text{d}s\\
    &~~= y_1^{\top}\hat{P}x_1+\myint x_1^{\top}\hat{Q}_1(s)y_2(s)\text{d}s+\myint y_2^{\top}(s) \hat{Q}_2(s) x_1 \text{d}s \\
    &~~~+ \myint y_2^{\top}(s) \hat{S}(s) x_2(s) \text{d}s + \myint\myinta{s}y_2^{\top}(s) \hat{R}_1(s,\eta)x_2(\eta) \text{d}\eta \text{d}s\\
    &~~~+\myint\myintb{s}y_2^{\top}(s) \hat{R}_2(s,\eta)x_2(\eta) \text{d}\eta \text{d}s\\
    &~~= \Big\langle{\opPpq{\hat{P}}{\hat{Q}}{\hat{R}}{\hat{S}}x},{y}\Big\rangle_{\R^{m}\times L_2^n[a, b]}\\
    &~~= \Big\langle{\opPpq{P}{Q}{R}{S}^*x},{y}\Big\rangle_{\R^{m}\times L_2^n[a, b]}.
\end{align*}}
where,
\begin{align*}
    &\hat{P} = P^{\top}, &\hat{S}(s) = S^{\top}(s), \\
    &\hat{Q}_1(s) = Q_2^{\top}(s), &\hat{R}_1(s,\eta) = R_2^{\top}(\eta,s), \\
    &\hat{Q}_2(s) = Q_1^{\top}(s), &\hat{R}_2(s,\eta) = R_1^{\top}(\eta,s). 
\end{align*}
This completes the proof. 
\QEDclosed
\end{pf}
\end{document}